\documentclass[12pt]{article}
\usepackage[cp1251]{inputenc}
\usepackage[russian]{babel}
\usepackage{amssymb,amsthm,amsmath}

\theoremstyle{plain}

\newtheorem{Proposition}{Утверждение}
\newtheorem{Theorem}{Теорема}
\newtheorem{THEOREM}{Теорема}

\newtheorem*{Corollary}{Следствие}
\newtheorem{Definition}{Определение}

\renewcommand\Im{\mathrm{Im\,}}
\renewcommand{\le}{\leqslant}
\renewcommand{\ge}{\geqslant}

\newcommand\wt{\widetilde}

\def\al{\alpha}
\def\la{\lambda}

\def\L{{\mathcal L}}

\def\bN{{\mathbb N}}
\def\bZ{{\mathbb Z}}
\def\bC{{\mathbb C}}
\def\bH{{\mathbb H}}

\def\P{\mathcal P}

\def\z{{\mathbf z}}
\def\y{{\mathbf y}}

\def\f{{\mathbf f}}

\def\R{{\mathfrak R}}

\newlength{\lenun}
\newlength{\lendu}
\begin{document}

{\bf УДК 517.984.52}

{\Large {\bf \centerline{$L_\mu\to L_\nu$ равносходимость спектральных разложений}} {\bf \centerline{ для системы
Дирака с $L_\varkappa$ потенциалом}}}

\medskip

{\bf \centerline{И.~В.~Садовничая\footnote{Работа выполнена при поддержке гранта РНФ 14-01-00754}}}
\bigskip

\medskip

В работе рассматривается оператор Дирака, порожденный в пространстве $\bH=L_2[0,\pi]\oplus L_2[0,\pi]\ni \y$
дифференциальным выражением
\begin{gather}\label{eq:lP}
\ell_P(\y)=B\y'+P\y,\quad \text{где}\\
 B = \begin{pmatrix} -i & 0 \\ 0 & i
        \end{pmatrix},
        \qquad
    P(x) = \begin{pmatrix} p_1(x) & p_2(x) \\ p_3(x) & p_4(x)
        \end{pmatrix},
        \qquad
     \y(x)=\begin{pmatrix}y_1(x)\\ y_2(x)\end{pmatrix}.\notag
\end{gather}
Функции $p_j$, $j=1,2,3,4$, предполагаются суммируемыми на отрезке $[0,\pi]$ и комплекснозначными. Общий вид краевых
условий для оператора $\L_{P,U}$ задается системой двух линейных уравнений $U(\y)=0$, где
\begin{equation*}
U(\y)=C\y(0)+D\y(\pi)=\begin{pmatrix}u_{11} & u_{12}\\ u_{21} &
u_{22}\end{pmatrix}\begin{pmatrix}y_1(0)\\
y_2(0)\end{pmatrix}+\begin{pmatrix}u_{13} & u_{14}\\ u_{23} &
u_{24}\end{pmatrix}\begin{pmatrix}y_1(\pi)\\
y_2(\pi)\end{pmatrix},
\end{equation*}
причем строки матрицы $\mathcal U:=(C,\,D)$ линейно независимы.

Система типа \eqref{eq:lP} была введена в рассмотрение П.~Дираком в 1929 году. Затем оператор Дирака изучался во многих
работах, но в основном  в случае  непрерывного потенциала  (см, например, монографию Б.~М.~Левитана и И.~С.~Саргсяна
\cite{LS} и литературу в ней).  Изучение асимптотики собственных значений оператора $\L_{P,U}$ опирается на
классические идеи и методы, восходящие к работам Г.~Биркгофа, Я.~Д.~Тамаркина и Р.~Лангера. В последнее время активно
изучаются операторы Дирака с негладкими потенциалами. Так, оператору $\L_P$ с периодическими и антипериодическими
условиями посвящена серия статей П.~Джакова и Б.~Митягина (см., например, работы \cite{DM11}, \cite{DM13} и литературу
в них). В статье А.~М.~Савчука и А.~А.~Шкаликова \cite{SavSh14} изучался случай $P\in L_\varkappa$,
$\varkappa\in[1,2]$. Случай негладкого потенциала рассматривался также в недавних работах М.~Ш.~Бурлуцкой,
В.~В.~Корнева, В.~П.~Курдюмова и А.~П.~Хромова (см., например, \cite{BKH1}).

В настоящей статье изучаются вопросы равносходимости спектральных разложений для возмущенного и невозмущенного
операторов Дирака. Тематика равносходимости восходит к работам У.~Дини, В.~А.~Стеклова, Я.~Д.~Тамаркина, А.~Хаара и
М.~Стоуна. В основном, изучалась равномерная на всем отрезке равносходимость. В работах В.~А.~Ильина и его
последователей (см. статьи \cite{Il}, \cite{Lom} и литературу в них) рассматривались вопросы равносходимости на
компактах внутри интервала $(0,1)$. В работе А.~М.~Гомилко и Г.~В.~Радзиевского \cite{Gom} получены оценки скорости
равносходимости для операторов Штурма--Лиувилля с классическими потенциалами. Подробный обзор классических результатов
о равносходимости содержится в монографии А.~Минкина \cite{Min}. Вопросы равносходимости для операторов Дирака с
негладкими потенциалами начали изучаться относительно недавно параллельно с изучением равносходимости для операторов
Штурма--Лиувилля с потенциалами--распределениями. По--видимому, первой работой здесь была статья В.~А.~Винокурова и
В.~А.~Садовничего \cite{VinSad}, в которой была доказана равномерная на всем отрезке $[0,\pi]$ равносходимость в
случае, когда раскладываемая функция $f\in L_1[0,\pi]$, а потенциал --- производная функции ограниченной вариации. В
работе автора \cite{Sad10} была доказана равномерная равносходимость в случае $f,\,u\in L_2[0,\pi]$ и краевых условий
Дирихле. П.~Джаков и Б.~Митягин в \cite{DM09} установили равномерную равносходимость для такого вида операторов
Штурма--Лиувилля с произвольными регулярными краевыми условиями, но для случая суммируемости в квадрате с
логарифмическим весом коэффициентов Фурье функции $f$. Ими также была доказана аналогичная теорема для системы Дирака с
потенциалом $P\in L_2[0,\pi]$.  В недавней статье А.~И.~Назарова, Д.~М.~Столярова и П.~Б.~Затицкого \cite{NSZ} был
установлен факт равносходимости для случая регулярного дифференциального оператора порядка $n\ge2$ с суммируемыми
коэффициентами и произвольными регулярными краевыми условиями. При этом раскладываемая функция $f\in L_1[0,\pi]$, а
равносходимость равномерна по шару $\|f\|_{L_1}\le R$.

Обозначим через \(J_{ij}\) определитель, составленный из \(i\)-го и \(j\)-го столбца матрицы $\mathcal U$.
\begin{Definition}\label{def:reg}
Краевые условия, определенные формой $U$, называются {\it регулярными} (по Биркгофу), если $J_{14}\cdot J_{23}\ne0$.
Оператор Дирака, порожденный регулярными краевыми условиями $U$, будем называть {\it регулярным}.
\end{Definition}
Рассмотрим вначале оператор $\L_{0,U}$, порожденный дифференциальным выражением $\ell_0(\y)=B\y'$ и регулярными
краевыми условиями.
\begin{THEOREM}[см., например, \cite{SavSh14}, теорема 3.1]  Спектр регулярного оператора $\L_{0,U}$ состоит из собственных
значений, которые можно записать двумя сериями $-\frac{i}{\pi}\ln z_0+2n$ и $-\frac{i}{\pi}\ln z_1+2n$, $n\in\mathbb
Z$, где $z_0$ и $z_1$ --- корни квадратного уравнения $ J_{23}z^2-[J_{12}+J_{34}]z-J_{14}=0$, а значения ветви
логарифма фиксируются в полосе $\Im\, z\in(-\pi,\pi]$. Можно занумеровать эти собственные значения одним индексом
$n\in\mathbb Z$: $\la^0_n=\begin{cases} \zeta_0+n,\quad \text{для четных }n,\\ \zeta_1+n,\quad \text{для нечетных
}n,\end{cases}\quad\text{где}\quad \zeta_0=-\frac{i}{\pi}\ln z_0,\quad \zeta_1=-\frac{i}{\pi}\ln z_1-1$.
\end{THEOREM}
Из теоремы i вытекает, что резольвента $\mathfrak{R}_0(\la)=(\L_{0,U}-\la I)^{-1}$ оператора $\L_{0,U}$ корректно
определена вне некоторой полосы $\Pi_a=\{\lambda\in\bC\,\vert\ |\Im \lambda|<a\}$ как оператор в пространстве $\bH$. В
\cite{SavSh14} (теорема 3.10) доказано также, что для некоторого числа $a>0$, зависящего только от краевых условий, вне
полосы $\Pi_a$ выполнена оценка: $\|\mathfrak{R}_0(\la)\|_{L_\mu\to L_\nu}\le C|\Im\la|^{-1+1/\mu-1/\nu}$. При этом
индексы $\mu$ и $\nu$ связаны соотношением $1\le\mu\le2\le\nu\le\infty$, а константа $C=C(U,\mu,\nu)$. Эту оценку можно
распространить на случай $1\le\mu\le\nu\le\infty$.
\begin{Proposition}\label{tm:R0est}
Пусть $\R_0(\la)$ --- резольвента регулярного оператора $\L_{0,U}$. Тогда для некоторого числа $a\ge1$, зависящего
только от краевых условий, вне полосы $\Pi_{a}=\{\la:|\Im\la|<a\}$ выполнена оценка
\begin{equation}\label{eq:R0est1}
\|\R_0(\la)\|_{L_\mu\to L_\nu}\le C_0|\Im\la|^{-1+1/\mu-1/\nu},\qquad\text{где}\ 1\le\mu\le\nu\le\infty.
\end{equation}
При этом число $C_0$ также зависит только от краевых условий.
\end{Proposition}
\begin{THEOREM} [\cite{DM11}, теорема 21] Пусть $\{\y_n^0\}_{n\in\bZ}$ --- система собственных и присоединенных функций
регулярного оператора $\L_{0,U}$, а $\{\z_n^0\}_{n\in\bZ}$ --- соответствующая биортогональная система. Для любой
непрерывной функции $\f$, имеющей ограниченную вариацию на $[0,\pi]$ и удовлетворяющей краевым условиям $U(\f)=0$, ряд
$\sum_{n=-\infty}^\infty\langle\f,\z_n^0\rangle\y_n^0$ сходится к $\f$ равномерно на отрезке $[0,\pi]$.
\end{THEOREM}
Заметим теперь, что случай потенциала произвольного вида \eqref{eq:lP} можно свести к случаю $p_1=p_4=0$. В дальнейшем,
чтобы не усложнять запись, мы будем писать $\f \in L_\al$, имея в виду, что $f_1\in L_\al[0,\pi]$ и $f_2\in
L_\al[0,\pi]$, или $P\in L_\al$, имея в виду, что все компоненты матрицы лежат в $L_\al[0,\pi]$.
\begin{Proposition}\label{st:potential} Пусть $P(x)$ --- произвольная матрица размера $2\times2$
с элементами $p_j\in L_1[0,\pi]$, $j=1,\,2,\,3,\,4$, а матрица $\mathcal{U}$ задает регулярные краевые условия. Тогда
оператор $\L_{P, U}$ подобен оператору $\L_{\wt P, \wt U}+\gamma I$, где
\begin{gather*}
\wt P(x)=\begin{pmatrix}0& \wt p_2(x)\\ \wt p_3(x)&0
\end{pmatrix},\notag\\
\wt p_2(x)= p_2(x)e^{i(\varphi(x)-\psi(x))},\qquad
\wt p_3(x)=p_3(x)e^{i(\psi(x)-\varphi(x))},\label{sim}\\
\varphi(x)=\gamma x-\int_0^x p_1(t)dt,\, \psi(x)=\int_0^x p_4(t)dt-\gamma x,\,
\gamma=\frac1{2\pi}\int_0^\pi( p_1(t)+ p_4(t))dt,\notag\\
\wt{\mathcal{U}}=(\wt C,\, \wt D),\quad  \wt C=C,\quad \wt D=\exp\left(\frac{i}2\int_0^\pi(p_4(t)- p_1(t))dt\right)
D.\notag
\end{gather*}
При этом если функция $P(\cdot)$ лежит в пространстве $L_\varkappa[0,\pi]$ для некоторого $\varkappa\in[1,\infty]$, то
и $\wt P(\cdot)\in L_\varkappa[0,\pi]$.
\end{Proposition}
Известно, что подобные операторы обладают одинаковыми спектральными свойствами: их спектры совпадают, системы корневых
функций образуют (или не образуют) базисы Рисса одновременно. Для нас важно также то, что операторы, осуществляющие
подобие, будут ограничены не только в пространстве $\bH$, но и в любом из пространств $L_\nu$, $\nu\ge 1$. Значит,
переход к подобным операторам не нарушает факта равносходимости. Однако, в случае потенциала общего вида \eqref{eq:lP}
потенциал невозмущенного оператора имеет, вообще говоря, диагональный вид.

Из теоремы ii вытекает, что система $\{\y_n^0\}_{n\in\bZ}$ собственных и присоединенных функций любого регулярного
оператора $\L_{0,U}$ полна в каждом из пространств $(L_\al[0,\pi])^2$, $\al\in[1,\infty)$. Можно доказать, что этим
свойством обладает и система корневых функций оператора $\L_{P,U}$ с потенциалом вида
\begin{equation}\label{eq:potential}
P(x)=\begin{pmatrix}0&p_2(x)\\p_3(x)&0\end{pmatrix}.
\end{equation}
\begin{Proposition}\label{tm:compl}
Пусть $\L_{P,U}$ --- регулярный оператор Дирака с суммируемым потенциалом $P$ вида \eqref{eq:potential}. Система
$\{\y_n\}_{n\in\bZ}$ его собственных и присоединенных функций полна в пространстве $(L_\al[0,\pi])^2$ для любого
$\al\in[1,\infty)$.
\end{Proposition}
Пусть $\{\la_n\}_{n\in\mathbb Z}$ --- система собственных значений регулярного оператора $\L_{P,U}$ с суммируемым
потенциалом вида \eqref{eq:potential}. Справедлива
\begin{THEOREM} [\cite{SavSad15}, теорема 4] Резольвента $\mathfrak{R}(\la)=(\L_{P, U}-\la
I)^{-1}$ любого регулярного оператора $\L_{P,U}$ с суммируемым комплекснозначным потенциалом $P$ вида
\eqref{eq:potential} определена при всех $\la\in\bC\setminus \{\la_n\}_{n\in\bZ}$ и является интегральным оператором в
$\bH$
\begin{equation*}
\mathfrak{R}(\la)\f=\int_0^\pi G(t,x,\la)\f(t)dt.
\end{equation*}
Функция $G(t,x,\la)=(g_{jk}(t,x,\la))$, $1\le j,\,k\le2$, непрерывна на квадрате $(t,x)\in[0,\pi]^2$ за исключением
диагонали.  Если считать потенциал $P$, краевые условия $U$ и число $\delta>0$ фиксированными, то вне семейства кружков
$\{|\la-\la_n|<\delta\}$ функция $G(t,x,\la)$ ограничена по модулю единой константой $M$.
\end{THEOREM}
Функцию $G(t,x,\la)$ мы будем называть функцией Грина оператора $\L_{P,U}$. В частности, функция Грина $G_0(t,x,\la)$
оператора $\L_{0,U}$ имеет вид
\begin{multline}\label{eq:Green0}
G_0(t,x,\la)=i\left(\frac{J_{12}}{\Delta_0(\lambda)}-\chi_{t>x}\right)\begin{pmatrix}e^{i\lambda(x-t)}&0\\0&-e^{i\lambda(t-x)}\end{pmatrix}+\\
+\frac{i}{\Delta_0(\lambda)}\begin{pmatrix}e^{i\lambda(x-\pi)}&0\\0&-e^{i\lambda(\pi-x)}\end{pmatrix}
\begin{pmatrix}J_{14}&J_{24}\\J_{13}&J_{23}\end{pmatrix}\begin{pmatrix}e^{-i\lambda t}&0\\0&-e^{i\lambda
t}\end{pmatrix},
\end{multline}
где через $\Delta_0(\la)$ обозначен характеристический определитель оператора $\L_{0,U}$.

В \cite{SavSh14} (теорема 4.3) доказано, что $\la_n=\la_n^0+o(1)$, $|n|\to\infty$. Воспользовавшись этим фактом и
теоремой i, можно построить семейство контуров $\gamma_n$ (здесь $|n|>N_0$ для некоторого $N_0\in\mathbb N$) такое, что
внутри каждого контура $\gamma_n$ лежат ровно два собственных значения оператора $\L_{P,U}$ ($\la_{2n}$ и
$\la_{2n+1}$). Обозначим
\begin{equation*}
\mathcal P_n:=\frac1{2\pi i}\int_{\gamma_n}\mathfrak{R}(\la)d\la,\qquad |n|>N_0, \qquad \mathfrak{R}(\la)=(\L_{P,U}-\la
I)^{-1}.
\end{equation*}
$\mathcal P_n$ является спектральным проектором на корневое подпространство, отвечающее собственным значениям
$\la_{2n}$ и $\la_{2n+1}$. Также обозначим через $\P_n^0$ спектральные проекторы на корневые подпространства оператора
$\L_{0,U}$, отвечающие собственным значениям $\la_{2n}^0$ и $\la_{2n+1}^0$. Построим контур $\Gamma_{m}$ так, чтобы он
охватывал в точности собственные значения $\la_n$ и $\la_n^0$ с номерами $n\in[-2m,2m+1]$. Положим $S_{m}:=\frac1{2\pi
i}\int_{\Gamma_{m}}\mathfrak R(\la)\,d\la,\quad m=N_0,\,N_0+1,\dots$. Из определения следует, что
$S_m=S_{N_0}+\sum_{|n|=N_0+1}^m\P_n$,  $m> N_0$. Заметим, что операторы $S_m$ можно задать в виде
\begin{equation*}
S_m(\f)=\sum_{|n|\le m}\Big[\langle\f,\z_{2n}\rangle\y_{2n}+\langle\f,\z_{2n+1}\rangle\y_{2n+1}\Big],
\end{equation*}
где $\{\y_n\}_{n\in\bZ}$ --- система собственных и присоединенных функций оператора $\L_{P,U}$, а $\{\z_n\}_{n\in\bZ}$
--- биортогональная система. Положим также
$$
S^0_m(\f)=\sum_{|n|\le m}\Big[\langle\f,\z^0_{2n}\rangle\y^0_{2n}+\langle\f,\z^0_{2n+1}\rangle\y^0_{2n+1}\Big],\quad
m=0,1,\dots.
$$
Легко видеть, что $S_m(\f)$ и $S_m^0(\f)$ представляют собой частичные суммы разложений функции $\f$ по системам
$\{\y_n\}_{n\in\bZ}$ и $\{\y_n^0\}_{n\in\bZ}$ соответственно. Основным результатом работы является
\begin{Theorem}
Пусть $\L_{P,U}$ --- регулярный оператор Дирака с потенциалом $P\in L_\varkappa[0,\pi]$ вида \eqref{eq:potential},
$\varkappa\in(1,\infty]$. Пусть $\f\in L_\mu[0,\pi]$, где $\mu\in[1,\infty]$. Тогда
\begin{equation*}
\|S_m(\f)-S_m^0(\f)\|_{L_{\nu}[0,\pi]}\longrightarrow0\quad \text{при}\ m\to+\infty,
\end{equation*}
если индекс $\nu\in[1,\infty]$ удовлетворяет неравенству
\begin{equation*}
\frac1\varkappa+\frac1\mu-\frac1\nu\le1
\end{equation*}
за исключением случая $\varkappa=\nu=\infty$, $\mu=1$.
\end{Theorem}
\begin{Corollary}
В случае оператора Дирака \eqref{eq:lP} с потенциалом $P$ общего вида и регулярными условиями $U=(C,\,D)$ утверждение
основной теоремы сохраняется, с той лишь разницей, что в определении последовательности $\{S_m^0\}$ в качестве системы
$\{\y_n^0\}$ следует взять систему собственных и присоединенных функций оператора $\L_{P_0,\wt U}$ с потенциалом
$$
P_0(x)=\begin{pmatrix}p_1(x)&0\\0&p_4(x)\end{pmatrix}
$$
и краевыми условиями
$$
\widetilde{\mathcal{U}}=(\wt C,\,\wt D),\quad \wt C=C,\quad\wt
D=\exp\left(\frac{i}2\int_0^\pi(p_4(t)-p_1(t))dt\right)D.
$$
\end{Corollary}

Схема доказательства теоремы заключается в следующем. Необходимо доказать сильную сходимость к нулю последовательности
операторов $\{S_m-S^0_m\}_{m\ge N_0}$, действующих из пространства $L_\mu[0,\pi]$ в $L_\nu[0,\pi]$. Для этого
достаточно проверить ограниченность норм этих операторов величиной, не зависящей от $m$, и сходимость
$\|S_m(\f)-S_m^0(\f)\|_{L_\nu}\to0$ на всюду плотном в $L_\mu$ множестве функций $\f$. В качестве этого множества мы
выберем линейную оболочку системы $\{\y_n\}_{n\in\bZ}$ собственных и присоединенных функций оператора $\L_{P,U}$ ---
для этого мы воспользуемся теоремой iii и утверждением \ref{tm:compl}. Наиболее сложным является доказательство
равномерной ограниченности семейства операторов $S_m-S_m^0$, $m\in\bN$. Здесь понадобится представление
$$
S_{m}-S^0_m:=\frac1{2\pi i}\int_{\Gamma_{m,b}} (\mathfrak R(\la)-\mathfrak R_0(\la))\,d\la,\qquad m=N_0,\,N_0+1,\dots
$$
и разложение
$$
\mathfrak R(\la)-\mathfrak R_0(\la)=-\R_0P\R_0+\R_0P\R_0P\R_0-\dots.
$$
В этом разложении, в свою очередь, наибольшие трудности составляет оценка первого слагаемого (остальные слагаемые
оцениваются единым образом с использованием неравенств \eqref{eq:R0est1}). Ключевую роль здесь играет явный вид функции
Грина невозмущенного оператора \eqref{eq:Green0}, оценки операторов типа Харди, а также теорема об интерполяции
линейных операторов.


\begin{thebibliography}{99}

\bibitem{LS} Левитан~Б.~М., Саргсян~И.~С. Операторы Штурма-Лиувилля и Дирака. М.: Наука,
1988.

\bibitem{DM11} Djakov~P., Mityagin~B. // Indiana Univ. Math. J.
2012. V. 61. \No 1. P. 359-398.

\bibitem{DM13} Djakov~P., Mityagin~B. // Proc. Amer. Math. Soc. 2013. V. 141. P. 1361-1375.

\bibitem{SavSh14} Savchuk~A.~M.,  Shkalikov~A.~A. // Math. Notes. 2014. V. 96 \No 5. P. 3-36.

\bibitem{BKH1} Бурлуцкая~М.~Ш., Курдюмов~В.~П., Хромов~А.~П. // Докл. Акад. Наук. 2012. Т. 443 \No 4. С. 414-417.

\bibitem{Il} Ильин~В.~А. // Дифф. уравнения. 1991. Т. 27 \No 4. С. 577-597.

\bibitem{Lom} Ломов~И.~С.  // Дифф. уравнения. 2001. Т. 37. \No 3. С. 328–342; \No 5. С. 648–660.

\bibitem{Gom} Гомилко~А.~М., Радзиевский~Г.~В. // Докл. РАН. 1991. Т. 316 \No 2. С. 265-270.

\bibitem{Min} Minkin~A. // J. Math. Sci. (New York). 1999. V. 96. P. 3631-3715.

\bibitem{VinSad} Винокуров~В.~А., Садовничий~В.~А. // Докл. РАН. 2001. Т. 380 \No 6. С. 731-735.

\bibitem{Sad10} Садовничая~И.~В. // Матем. сб. 2010. Т. 201 \No 9. С. 61-76.

\bibitem{DM09} Djakov~P., Mityagin~B. // Contemp. Math. 2009. V. 481. P. 59-80.

\bibitem{NSZ} Nazarov~A.~I., Stolyarov~D.~M., Zatitskiy~P.~B. // J. of Spectral Th. 2014. V. 4 \No 2. P. 365-389.

\bibitem{SavSad15} Савчук~А.~М., Садовничая~И.~В. //
Современная математика. Фундаментальные направления. 2015. Т. 58.

\end{thebibliography}
\end{document}